\documentclass[a4paper, 11pt,twoside]{article}

\usepackage{color}
\usepackage{amssymb}
\usepackage{amsmath}
\usepackage{bbm}
\usepackage{mathrsfs}
\usepackage{tikz}
\usepackage{caption}
\usepackage{amsthm}
\usetikzlibrary {positioning}

\newtheorem{theorem}{Theorem}
\newtheorem{corollary}{Corollary}
\newtheorem{lemma}{Lemma}
\newtheorem{conjecture}{Conjecture}
\newtheorem{example}{Example}
\newtheorem{proposition}{Proposition}
\newtheorem{property}{Property}

\newtheorem{definition}{Definition}
\begin{document}

\author{\footnotesize\textbf{Henk Don}\\
\\
\footnotesize Institute for Mathematics, Astrophysics, and Particle Physics\\
\footnotesize Faculty of Science, Radboud University Nijmegen;
e-mail:henkdon@gmail.com}

\title{The \v{C}ern\'y conjecture and $1-$contracting automata}
\maketitle

\begin{abstract}
A deterministic finite automaton is synchronizing if there exists
a word that sends all states of the automaton to the same state.
\v{C}ern\'y conjectured in 1964 that a synchronizing automaton
with $n$ states has a synchronizing word of length at most
$(n-1)^2$. We introduce the notion of aperiodically
$1-$contracting automata and prove that in these automata all
subsets of the state set are reachable, so that in particular they
are synchronizing. Furthermore, we give a sufficient condition
under which the \v{C}ern\'y conjecture holds for aperiodically
$1-$contracting automata. As a special case, we prove some results
for circular automata.
\end{abstract}

\section{Introduction}

Let $\mathscr{A} = (Q,\Sigma,\delta)$ be a deterministic finite
automaton (DFA), where $Q$ denotes the state set, $\Sigma$ the input
alphabet, and $\delta:Q\times \Sigma\rightarrow Q$ the transition
function. We denote the set of finite words over $\Sigma$ by
$\Sigma^\star$. The transition function $\delta$ extends uniquely
to a function $\delta: Q\times \Sigma^\star\rightarrow Q$.

The automaton $\mathscr{A}$ is called \emph{synchronizing} if
there exists a word $w\in\Sigma^\star$ and $q\in Q$ such that $\delta(q',w) =
q$ for all $q'\in Q$. The word $w$ is then said to be a \emph{synchronizing
word} for $\mathscr{A}$.

The following longstanding conjecture is due to \v{C}ern\'y
(\cite{cerny}, 1964):
\begin{conjecture} If $\mathscr{A}$ is a synchronizing $n$-state automaton,
then there exists a synchronizing word for $\mathscr{A}$ of length
at most $(n-1)^2$.
\end{conjecture}

\v{C}ern\'y constructed for every $n$ an $n$-state automaton
$\mathscr{C}_n$ in which the bound of his conjecture is attained.
We show this automaton in the next example for the case $n=4$.

\begin{example}  \emph{Consider the automaton $\mathscr{C}_4$ defined by the transition graph shown in Figure \ref{fig:transgraph1}. The word $w = baaabaaab$ (having length $9$) is synchronizing, since it maps all states to state $1$. Moreover, $w$ is the shortest synchronizing
word for $\mathscr{C}_4$.}\hfill $\blacksquare$

\begin{center}
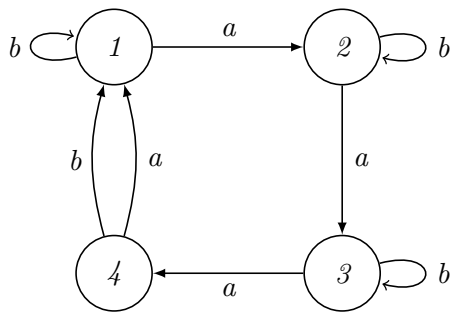

\begin{tikzpicture}[-latex ,auto ,node distance =3 cm and 3cm ,on grid ,
semithick ,
state/.style ={ circle ,top color =white , bottom color = white!20 ,
draw, black , text=black , minimum width =1 cm}]
\node[state] (1) {1};
\node[state] (2) [right =of 1] {2};
\node[state] (3) [below =of 2] {3};
\node[state] (4) [left  =of 3] {4};
\path (1) edge [bend right = 0] node[above] {a} (2);
\path (1) edge [loop left] node[left] {b} (1);
\path (2) edge [bend right = 0] node[right] {a} (3);
\path (2) edge [loop right] node[right] {b} (2);
\path (3) edge [bend right = 0] node[below] {a} (4);
\path (3) edge [loop right] node[right] {b} (3);
\path (4) edge [bend left = -15] node[right] {a} (1);
\path (4) edge [bend right = -15] node[left] {b} (1);
\end{tikzpicture}\captionof{figure}{Transition graph of
$\mathscr{C}_4$.}\label{fig:transgraph1}
\end{center}

\end{example}

\v{C}ern\'y himself provided an upper bound of $2^n-n-1$ for the
length of the shortest synchronizing word. A substantial
improvement was given by Starke \cite{starke}, who was the first
to give a polynomial upper bound, namely
$1+\frac{1}{2}n(n-1)(n-2)$. The best known upper bound is still a
cubic one, $\frac{1}{6}(n^3-n)$, that was established by Pin in
1983 \cite{pin}. He reduced proving this upper bound to a purely
combinatorial problem which was then solved by Frankl \cite{frankl}.
Since then for more than 30 years no progress for the general case
has been made.

The conjecture has been proved for some particular classes of
automata, such as circular automata, one-cluster automata and
aperiodic automata. For these results and some more partial
answers, see \cite{dubuc,beal,trahtman,eppstein, volkov2, kari,
almeida}. For a survey on synchronizing automata and the
\v{C}ern\'y conjecture, we refer to \cite{volkov}.

In this paper we look at $n-$state automata in which every $(n-1)-$subset of the state set $Q$ is reachable from $Q$. Such automata will be called $1-$contracting. A word with the property that it maps $Q$ to an $(n-1)-$subset of $Q$ is called $1-$deficient. If $w$ is a $1-$deficient word, the state that is not in the image of $w$ is said to be the \emph{excluded} state. There also must be a unique state in the image which is reached twice by $w$. This state will be called the \emph{contracting} state for $w$. 

In a $1-$contracting automaton, for every state $q$ there exists a $1-$deficient word that excludes $q$. A collection $W$ of words is called $1-$contracting if for all $q$ it contains exactly one word $w_q$ which excludes $q$. To such a collection we can associate a function $\sigma_W$ on $Q$ that maps each state $q$ to the unique contracting state for $w_q$. This function will be called the state map induced by $W$. If for some $1-$contracting collection $W$ the state map is a cyclic permutation on $Q$, then the automaton is called aperiodically $1-$contracting. For formal definitions and
some more details, we refer to Section \ref{sec:hamiltonian}. 

Our main results are presented in Section \ref{sec:synchronization}, where we will prove that aperiodically $1-$contracting automata are synchronizing and that every nonempty subset of the state set is reachable. By imposing a (weak) restriction on the length of the $1-$deficient words, we show that each $k-$subset is reachable by a word of length at most $n(n-k)$ which implies that the \v{C}ern\'y conjecture holds true in this case.

In Section \ref{section:circular}, we give some examples and discuss the scope of our results and their relation to circular automata. For circular automata it was shown by Dubuc (\cite{dubuc}) that if they are synchronizing, then the shortest synchronizing word has length at most $(n-1)^2$. We give sufficient conditions for circular automata to be aperiodically $1-$contracting. This means that for a subclass of circular automata, we prove synchronization. For this subclass, we also prove that the \v{C}ern\'y conjecture is satisfied and moreover that also here each $k-$subset is reachable by a word of length at most $n(n-k)$.

The class of aperiodically $1-$contracting automata also contains automata satisfying a weaker notion of circularity: automata with a circle that uses multiple labels and goes through all states. We give an example that illustrates this and note that there is an infinite subfamily of such automata in which the \v{C}ern\'y conjecture holds. Finally, we provide an example of an automaton that has no circle through all states, but still satisfies the \v{C}ern\'y conjecture because it falls into the framework of aperiodically $1-$contracting automata, which demonstrates the extensiveness of this class. We finish our paper with some open questions and a conjecture.

\section{Notation and preliminaries}

\subsection{The power automaton}\label{sec:power}

The power automaton $\mathcal{P}(\mathscr{A})$ of a DFA
$\mathscr{A} = (Q,\Sigma,\delta)$ has state set $2^Q$ and its
transition function (with a slight abuse of notation again denoted
by $\delta$) is uniquely determined by the transition function of
$\mathscr{A}$:
$$
\delta(S,w) = \left\{\delta(q,w)\mid q\in S\right\}.
$$
The automaton $\mathscr{A}$ is synchronizing if and only if there
is a path from $Q$ to a singleton in the transition graph of
$\mathcal{P}(\mathscr{A})$. For some
$S\subseteq Q$, if $\delta(Q,w) = S$, we say that $S$ is \emph{reachable} by $w$ (from $Q$).

\begin{example}\label{example:power} \emph{The power automaton of $\mathscr{C}_4$ is given in Figure \ref{fig:power}, omitting selfloops. The synchronizing word $w = baaabaaab$ in $\mathscr{C}_4$ corresponds to the
double arrowed path. Note that from the triple
$\left\{1,3,4\right\}$ we can already reach a pair. So we could be
more greedy and try to synchronize with a word that starts as
$\tilde w = baab\ldots$. However, this would lead to the
synchronizing word $\tilde w = baababaaab$, which in fact is
\emph{longer} than $w$.}\hfill$\blacksquare$
\end{example}

\begin{center}
\begin{tikzpicture}[-latex ,auto ,node distance =1.65 cm and 1.65cm ,on grid ,
semithick , state/.style ={ circle ,top color =white , bottom
color = white!20 , draw, black , text=black , minimum width =1
cm}]

\node[state] (234) {\underline{2}34};%
\node[state] (123) [above right =of 234] {\underline{1}23};%
\node[state] (124) [below right =of 123] {12\underline{4}};%
\node[state] (12) [right =of 124] {12};%
\node[state] (14) [right =of 12] {14};%
\node[state] (1) [right =of 14] {1};%
\node[state] (4) [above right =of 1] {4};%
\node[state] (3) [below right =of 4] {3};%
\node[state] (1234) [above =of 12] {Q=1234};%
\node[state] (134) [below right =of 234] {1\underline{3}4};%
\node[state] (2) [below right =of 1] {2};%
\node[state] (23) [below =of 12] {23};%
\node[state] (34) [below =of 14] {34};%
\node[state] (24) [right =of 134] {24};%
\node[state] (13) [below =of 24] {13};%

\path (1234) edge [double distance =.5pt] node[above] {b} (123);%
\path (123) edge [double distance =.5pt] node[below right] {a} (234);%
\path (234) edge [double distance =.5pt] node[above right] {a} (134);%
\path (134) edge [double distance =.5pt] node[above left] {a} (124);%
\path (124) edge [bend right = 0] node[below left] {a} (123);%

\path (12) edge [double distance =.5pt] node[right] {a} (23);%
\path (23) edge [double distance =.5pt] node[above] {a} (34);%
\path (34) edge [double distance =.5pt] node[left] {a} (14);%
\path (14) edge [bend right = 0] node[below] {a} (12);%

\path (1) edge [bend right = 0] node[above right] {a} (2);%
\path (2) edge [bend right = 0] node[above left] {a} (3);%
\path (3) edge [bend right = 0] node[below left] {a} (4);%
\path (4) edge [bend right = 0] node[below right] {a} (1);%

\path (24) edge node[right] {a} (13);%
\path (13) edge (24);%

\path (234) edge [bend left = 15] node[above left] {b} (123);%
\path (24) edge [bend left = 0] node[above left] {b} (12);%
\path (124) edge [double distance =.5pt] node[above] {b} (12);%
\path (4) edge [bend right = 15] node[above left] {b} (1);%
\path (14) edge [double distance =.5pt] node[above] {b} (1);%
\path (134) edge [bend left = 0] node[below left] {b} (13);%
\path (34) edge [bend left = 0] node[below] {b} (13);%


\end{tikzpicture}
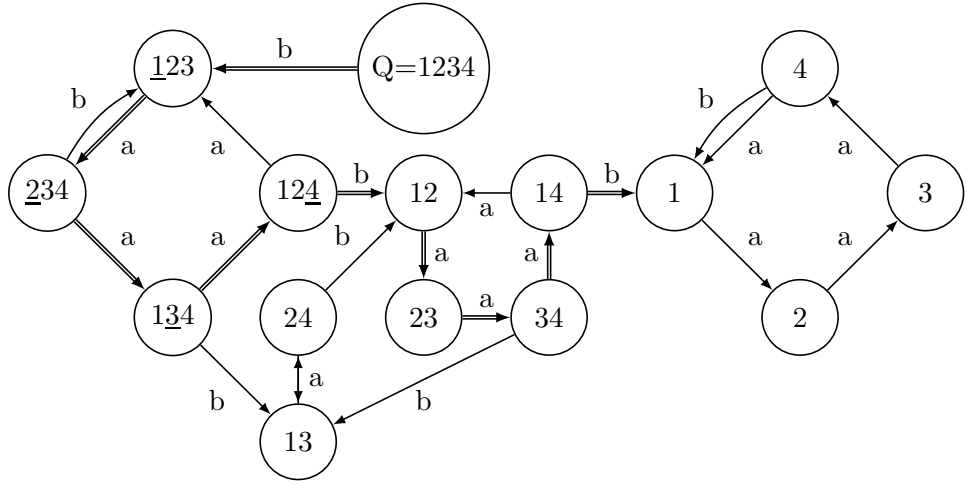
\captionof{figure}{Transition graph of
$\mathcal{P}(\mathscr{C}_4$). Double arrows indicate the shortest
path from $Q$ to a singleton. Unique contracting states on this path are underlined.}\label{fig:power}
\end{center}

\subsection{Some remarks on the transition function}

For the empty word $\varepsilon$, we define the transition
function to be the identity: $\delta(q,\varepsilon) = q$. For
$S\subseteq Q$, we define the inverse of the transition function
under some word $w$ by
$$
\delta^{-1}(S,w) = \left\{q \mid \delta(q,w)\in S\right\}.
$$
For $w\in\Sigma^\star$ and $q\in Q$, we define the $w-$indegree of
$q$ to be $|\delta^{-1}(q,w)|$. If $|\delta^{-1}(q,w)|\geq 2$, we
say that $q$ is a \emph{contracting state} for $w$. Some useful
properties of the transition function and its inverse are listed
here:
\begin{property}\label{prop:transfunc} Let $\mathscr{A} = (Q,\Sigma,\delta)$ be a DFA. The (extended) transition function $\delta:2^Q \times \Sigma^\star \rightarrow 2^Q$ and its inverse satisfy the following:
\begin{enumerate}
\item For all non-empty $S\subseteq Q$ and $w\in\Sigma^\star$,
$$
1\leq |\delta(S,w)| \leq |S|.
$$
\item For all $w\in\Sigma^\star$,
$$
\sum_{q\in Q}|\delta^{-1}(q,w)| = |Q|.
$$
\item For all $S_1\subseteq S_2 \subseteq Q$ and $w\in\Sigma^\star$,
$$
|S_2| - |\delta(S_2,w)| \geq |S_1| - |\delta(S_1,w)| \geq 0.
$$
\item For all $S_1, S_2\subseteq Q$ and $w\in\Sigma^\star$,
$$
\delta(S_1,w)\cup \delta(S_2,w) = \delta(S_1\cup S_2,w).
$$
\end{enumerate}
\end{property}

\subsection{$1-$contracting automata}\label{sec:hamiltonian}

In this paper, we study automata in which each $(n-1)-$subset of
the state set is reachable from $Q$. In this section we formally introduce
the class of automata we will consider and give some related
terminology.

\begin{definition}
Let $\mathscr{A} = (Q,\Sigma,\delta)$ be a DFA. If a word
$w_q\in\Sigma^\star$ satisfies
$$
\delta(Q,w_q) = Q\setminus\left\{q\right\},
$$
then $w_q$ will be called a $1-$deficient word that excludes $q$.
\end{definition}

\begin{definition}
Let $\mathscr{A} = (Q,\Sigma,\delta)$ be a DFA with $n$ states. A
collection $W\subseteq\Sigma^\star$ will be called a
$1-$contracting collection for $\mathscr{A}$ if for all $q\in Q$
it contains exactly one $1-$deficient word that excludes $q$. If
such a collection exists, $\mathscr{A}$ is called $1-$contracting.
If $|w|\leq n$ for all $w\in W$, then $W$ is called an efficient
$1-$contracting collection.
\end{definition}

\begin{lemma}\label{lemma:efficient} Every $1-$contracting automaton $\mathscr{A} =
(Q,\Sigma,\delta)$ admits an efficient $1-$contracting collection.
\end{lemma}

\textbf{Proof.} Suppose $w_q$ is a $1-$deficient word that
excludes $q$. This word defines a path in the power automaton
$\mathcal{P}(\mathscr{A})$ that visits only sets of size at least
$n-1$. After erasing loops, this path has length at most $n$,
which means that $w_q$ can be reduced to have length at most $n$.
\hfill $\Box$\\
\\
The previous lemma states that in a $1-$contracting automaton on
$n$ states, every $(n-1)-$set is reachable by a word of length at
most $n$. In such automata, also many $(n-2)$-sets (in fact at
least $\binom{n}{2}-n$) are reachable by words of length at most
$2n$. The main goal is to show inductively that this behavior is
inherited by sets of arbitrary size, i.e. that many $(n-k)-$sets
are reachable by words of length at most $kn$. It turns out that
it is not sufficient to merely have $1-$contractivity, we also
need some kind of aperiodicity, as the following example shows.

\begin{example}\label{example:fail}
\emph{Let $\mathscr{A}$ be the automaton with the transition graph
given in Figure \ref{fig:transgraph2}. The collection $W =
\left\{c,ca,cab,cb\right\}$ is an efficient $1-$contracting
collection for $\mathscr{A}$ since
\begin{eqnarray*}
\delta(Q,c) = Q\setminus\left\{1\right\},&\qquad & \delta(Q,cab) = Q\setminus\left\{3\right\},\\
\delta(Q,ca) = Q\setminus\left\{2\right\},&\qquad & \delta(Q,cb) =
Q\setminus\left\{4\right\}.
\end{eqnarray*}
Thus, $\mathscr{A}$ is $1-$contracting. However, $\mathscr{A}$ has
no synchronizing word. For instance there does not exist a word
$w\in\Sigma^\star$ such that $\delta(\left\{1,2\right\},w)$ is a
singleton.}\hfill$\blacksquare$
\begin{center}
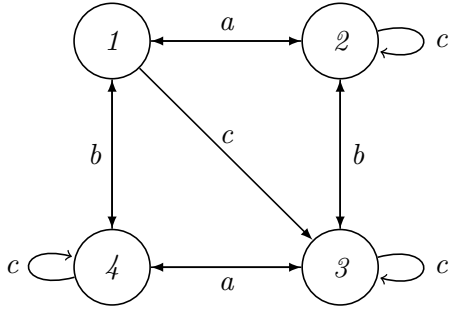

\begin{tikzpicture}[-latex ,auto ,node distance =3 cm and 3cm ,on grid ,
semithick , state/.style ={ circle ,top color =white , bottom
color = white!20 , draw, black , text=black , minimum width =1
cm}]

\node[state] (1) {1}; %
\node[state] (2) [right =of 1] {2};
\node[state] (3) [below =of 2] {3}; %
\node[state] (4) [left  =of 3] {4};

\path (1) edge [bend right = 0] node[above] {a} (2); \path (2)
edge [bend right = 0] node[above] {} (1); \path (3) edge [bend
right = 0] node[below] {a} (4); \path (4) edge [bend right = 0]
node[below] {} (3);
\path (2) edge [bend right = 0] node[right] {b} (3); \path (3)
edge [bend right = 0] node[below] {} (2); \path (1) edge [bend
right = 0] node[left] {b} (4); \path (4) edge [bend right = 0]
node[below] {} (1);
\path (1) edge [bend right = 0] node[above] {c} (3); \path (2)
edge [loop right] node[right] {c} (2); \path (3) edge [loop right]
node[right] {c} (3); \path (4) edge [loop left] node[left] {c}
(4);
\end{tikzpicture}
\captionof{figure}{A $1-$contracting automaton that fails to
synchronize.}\label{fig:transgraph2}
\end{center}
\end{example}

Our notion of aperiodicity will be formulated in terms of
indegrees corresponding to the $1-$deficient words. Suppose
$\mathscr{A}$ is a $1-$contracting automaton with $n$ states. Let
$q\in Q$ and let $w_q\in\Sigma^\star$ be a $1-$deficient word that
excludes $q$. Then $|\delta^{-1}(q,w_q)| = 0$ and
$|\delta^{-1}(q',w_q)| \geq 1$ for all $q'\neq q$. Since
$\sum_{q\in Q}|\delta^{-1}(q,w_q)| = n$, there is a unique
contracting state $q^c$ for $w_q$ which has $w_q-$indegree 2. See
also Figure \ref{fig:power} for an illustration.

So in an automaton $\mathscr{A}$, a $1-$deficient word $w$ has a
unique excluded state and a unique contracting state.

\begin{definition}
Let $\mathscr{A}$ be a DFA and let $W\subseteq \Sigma^\star$ be a
$1-$contracting collection for $\mathscr{A}$. Define the following
collection of pairs of states $\Pi_W\subseteq Q\times Q$:
$$
\Pi_W:=\left\{(q,q^c):\exists w\in W \textrm{\ s.t.\ }
|\delta^{-1}(q,w)| = 0, |\delta^{-1}(q^c,w)| = 2\right\}.
$$
Let the map $\sigma_W$ from $Q$ to $Q$ be defined by
$$
\sigma_W(q) = q^c \Leftrightarrow (q,q^c)\in \Pi_W.
$$
The map $\sigma_W$ will be called the state map induced by $W$.
\end{definition}
Since for every $q\in Q$ there is exactly one $w\in W$ which
excludes $q$, the state map induced by $W$ is well-defined.

\begin{example}\label{example:statemap}
\emph{Consider the automaton $\mathscr{A}$ of Figure
\ref{fig:transgraph2} with the $1-$contracting collection $W =
\left\{c,ca,cab,cb\right\}$. The contracting states for the words
in this collection follow directly from the equalities
$$
|\delta^{-1}(3,c)| = |\delta^{-1}(4,ca)| =  |\delta^{-1}(1,cab)| =
|\delta^{-1}(2,cb)| = 2.
$$
So in this case $\Pi_W = \left\{(1,3),(2,4),(3,1),(4,2)\right\}$,
and therefore the state map $\sigma_W:Q\rightarrow Q$ is defined
by
$$
\sigma_W(1) = 3,\qquad \sigma_W(2) = 4,\qquad \sigma_W(3) =
1,\qquad \sigma_W(4) = 2.
$$}\hfill$\blacksquare$
\end{example}

We now have the ingredients to formulate our aperiodicity
property.

\begin{definition}\label{def:aperiodic} Let $\mathscr{A} = (Q,\Sigma,\delta)$ be a $1-$contracting
automaton. If there exists a $1-$contracting collection
$W\subseteq \Sigma^\star$ for which the induced state map
$\sigma_W$ is a cyclic permutation on $Q$, then $\mathscr{A}$ is called
aperiodically $1-$contracting.
\end{definition}

If $Q$ has $n$ states, then a cyclic permutation $\sigma$ on $Q$
is characterized by the fact that $\sigma^n(q) = q$ for all $q\in
Q$ and that $\sigma^k(q)\neq q$ if $1\leq k < n$. So this is
independent of the way the states in $Q$ are enumerated. As can be
easily checked, the state map $\sigma_W$ of Example
\ref{example:statemap} is not a cyclic permutation.

\section{Synchronization of $1-$contracting automata}\label{sec:synchronization}

\subsection{Main results}

In this section we prove that in aperiodically $1-$contracting automata all subsets of the state set are reachable by concatenating $1-$deficient words. As a direct consequence, such automata are synchronizing.
\begin{theorem}\label{maintheorem1} Let $\mathscr{A} = (Q,\Sigma,\delta)$ be an aperiodically $1-$contracting DFA with $n$
states. Let $W = \left\{w_1,\ldots,w_n\right\}\subseteq
\Sigma^\star$ be a corresponding $1-$contracting collection for
which the state map $\sigma_W$ is a cyclic permutation on $Q$.

Then every nonempty $k-$subset $S$ of $Q$ is reachable by a word
$w_S$ of the form $w_S = w_{i_1}\ldots w_{i_{n-k}}$.
\end{theorem}
\begin{corollary} Under the conditions of Theorem
\ref{maintheorem1}, $\mathscr{A}$ is synchronizing and there exist
synchronizing words of the form $w = w_{i_1}\ldots w_{i_{n-1}}$.
\end{corollary}
Furthermore, we show that the \v{C}ern\'y conjecture holds true if
there is an \emph{efficient} $1-$contracting collection $W$ for
which the state map is a cyclic permutation. By Lemma
\ref{lemma:efficient}, for $W$ to be efficient, it is sufficient
if the $1-$deficient words can be chosen to have minimal length.

\begin{theorem}\label{maintheorem2}
Let $\mathscr{A} = (Q,\Sigma,\delta)$ be an aperiodically
$1-$contracting DFA with $n$ states. If there exists an efficient
$1-$contracting collection $W\subseteq \Sigma^\star$ for which
$\sigma_W$ is a cyclic permutation on $Q$, then
\begin{enumerate}
\item The shortest synchronizing word of $\mathscr{A}$ has length at most $(n-1)^2$.
\item For every nonempty set $S\subseteq Q$ of size $k$, there exists a word $w_S$ of length at most $n(n-k)$ such that $\delta(Q,w_S) = S$.
\end{enumerate}
\end{theorem}

\subsection{Proofs of main results}

Before turning to the proofs, we introduce some terminology. Let
$\mathscr{A}$ be an aperiodically $1-$contracting DFA. Denote the
state set $Q$ by $\left\{q_1,\ldots,q_n\right\}$ and choose a
corresponding $1-$contracting collection $W =
\left\{w_1,\ldots,w_n\right\}$ such that the state map $\sigma_W$
is a cyclic permutation on $Q$. Define $q_i^c$ to be the contracting
state for the word $w_i$ that excludes $q_i$, so $q_i^c =
\sigma_W(q_i)$. The strategy of the proof will be to show that
every $k-$subset of $Q$ can be reached by properly choosing $n-k$
of the words $w_i$ and concatenating them. The following
collections of subsets of $Q$ play an important role. We let
$\mathscr{S}_{0,i} =\mathscr{T}_0 =\mathscr{U}_{0,i}
=\left\{Q\right\}$ and for $k=1,\ldots ,n-1$, we recursively
define the collections $\mathscr{S}_k$, $\mathscr{T}_k$ and
$\mathscr{U}_k$:
\begin{eqnarray*}
\mathscr{S}_{k,i} &=& \left\{\delta(S,w_i)\mid
S\in\mathscr{U}_{k-1,i}\right\},\\
\mathscr{T}_k &=& \bigcup_{i=1}^n \mathscr{S}_{k,i},\\
\mathscr{U}_{k,i} &=& \left\{S\in\mathscr{T}_k \mid
\delta^{-1}(q_i^c,w_i)\subseteq S\right\}.
\end{eqnarray*}
The following sequence of results gives properties of these collections that together are sufficient to prove Theorems \ref{maintheorem1} and \ref{maintheorem2}.

\begin{lemma}\label{lemma:length}
Let $S\in \mathscr{T}_k$ for $1\leq k\leq n-1$. Then there exists a word $w_S\in\Sigma^\star$ of the form $w_S = w_{i_1}\ldots w_{i_{k}}$ such that $\delta(Q,w_S) = S$.
\end{lemma}
\textbf{Proof.} If $S\in \mathscr{S}_{1,i}$, then $S = \delta(Q,w_i) = Q\setminus\left\{q_i\right\}$, so we can take $w_S$ to be the $1-$deficient word $w_i$. So the statement holds for $k=1$.\\
\\
Assume that the statement holds for some $k\geq 1$ and take $S\in\mathscr{S}_{k+1,i}$. Then $\delta^{-1}(S,w_i)\in\mathscr{U}_{k,i}$, and therefore $\delta^{-1}(S,w_i)\in\mathscr{T}_{k}$. By definition of $\mathscr{T}_{k}$, there exists $1\leq j\leq n$ such that $\delta^{-1}(S,w_i)\in\mathscr{S}_{k,j}$.\\
\\
By the assumption, there exists a word $\tilde w$ of the form $\tilde w = w_{i_1}\ldots w_{i_{k}}$ such that $\delta(Q,\tilde w) = \delta^{-1}(S,w_i)$. Therefore, the word $w_S:=\tilde ww_i$ has the desired form and satisfies $\delta(Q,w_S) = S$.
\hfill $\Box$

\begin{lemma}\label{lemma:size} Let $S\subseteq Q$ and $0\leq k\leq n-1$. If $S\in\mathscr{S}_{k,i}$ for some $i$, $S\in\mathscr{U}_{k,i}$ for some $i$ or $S\in\mathscr{T}_k$, then $|S| = n-k$.
\end{lemma}
\textbf{Proof.} For $k=0$, the statement is obvious. Now take
$k_0>0$ and suppose the statement is true for $0\leq k
<k_0$.\\
\\
It suffices to prove for all $i$ that $|S|=n-k$ if $S\in\mathscr{S}_{k_0,i}$, the other statements then follow from the definitions of $\mathscr{T}_k$ and $\mathscr{U}_{k,i}$. So let $S\in\mathscr{S}_{k_0,i}$. Then there exists $S'\in\mathscr{U}_{k_0-1,i}$ such that $S = \delta(S',w_i)$. The
induction hypothesis gives $|S'| = n-k_0+1$. By definition of
$\mathscr{U}_{k_0-1,i}$ we have $\delta^{-1}(q_i^c,w_i)\subseteq S'$,
which implies $|S|\leq |S'|-1 = n-k_0$. The other inequality
follows from the third property of the transition function:
$$
1 = |Q| - |\delta(Q,w_i)| \geq |S'| - |\delta(S',w_i)| =
n-k_0+1-|S|,
$$
which implies $|S| \geq n-k_0$. Therefore $|S| = n-k_0$.\hfill $\Box$

\begin{lemma}\label{lemma:bijection}
For all $i$, the function $\delta(\cdot,w_i)$ defines a bijection from $\mathscr{U}_{k-1,i}$ to $\mathscr{S}_{k,i}$.
\end{lemma}
\textbf{Proof.} Surjectivity follows immediately from the definition of $\mathscr{S}_{k,i}$. To prove injectivity, take $S_1, S_2\in \mathscr{U}_{k-1,i}$ with $S_1\neq S_2$. By Lemma \ref{lemma:size}, $|S_1| = |S_2| = n-k+1$, so the size of their union satisfies $|S_1\cup S_2| \geq n-k+2$. Application of property three of the transition function gives
$$
|S_1\cup S_2|-|\delta(S_1\cup S_2,w_i)| \leq |Q|-|\delta(Q,w_i)| = 1.
$$
Therefore,
$$
|\delta(S_1,w_i)\cup \delta(S_2,w_i)| = |\delta(S_1\cup S_2,w_i)| \geq |S_1\cup S_2| - 1 \geq n-k+1.
$$
Lemma \ref{lemma:size} gives $|\delta(S_1,w_i)|= |\delta(S_2,w_i)| = n-k$, hence it follows that $\delta(S_1,w_i)\neq \delta(S_2,w_i)$.\hfill $\Box$\\
\\
The following proposition shows that the collections $\mathscr{S}_{k,i}$ and $\mathscr{T}_k$ are maximal in the sense that they contain every subset of $Q$ that is not immediately excluded by definition.

\begin{proposition}\label{prop:maxsets}
For $1\leq k\leq n-1$ and $1\leq i\leq n$, let
\begin{eqnarray*}
\overline{\mathscr{S}}_{k,i} &=& \left\{S\subseteq Q \mid q_i^c \in S, q_i\notin S, |S| = n-k\right\}.\\
\overline{\mathscr{T}}_k &=& \left\{S\subseteq Q \mid |S| = n-k\right\}.
\end{eqnarray*}
Then for all $i$ and $k$, the collections $\mathscr{S}_{k,i}$ and $\mathscr{T}_k$ satisfy
$$
\mathscr{S}_{k,i} = \overline{\mathscr{S}}_{k,i}\qquad\textrm{and}\qquad\mathscr{T}_k = \overline{\mathscr{T}}_k.
$$
\end{proposition}

\textbf{Proof.} Let $1\leq k\leq n-1$. For $S\in\mathscr{S}_{k,i}$, its definition gives that $q_i^c\in S$ and also
$$
S \subseteq \delta(Q,w_i) = Q\setminus\left\{q_i\right\}.
$$
Lemma \ref{lemma:size} gives $|S| = n-k$. Therefore, $S\in\overline{\mathscr{S}}_{k,i}$ and hence $\mathscr{S}_{k,i} \subseteq \overline{\mathscr{S}}_{k,i}$. Clearly we also have $\mathscr{T}_k \subseteq \overline{\mathscr{T}}_k$.\\
\\
For the reversed inclusions, we will give a proof by induction. First we check them for $k=1$. Using the definitions of $\mathscr{S}_{1,i}$, $\mathscr{U}_{0,i}$, we obtain
$$
\mathscr{S}_{1,i} = \left\{\delta(S,w_i)\mid
S\in\mathscr{U}_{0,i}\right\} = \left\{\delta(Q,w_i)\right\} = \left\{Q\setminus\left\{q_i\right\}\right\} = \overline{\mathscr{S}}_{1,i}.
$$
This also gives $\mathscr{T}_1 = \bigcup_{i=1}^n \left\{Q\setminus\left\{q_i\right\}\right\} = \overline{\mathscr{T}}_1$.\\
\\
Now let $k_0>1$ and suppose the two statements are true for $1\leq k<k_0$. Since $\mathscr{T}_{k_0-1} = \overline{\mathscr{T}}_{k_0-1}$, we know that
$$
\mathscr{U}_{k_0-1,i} = \left\{S\subseteq Q \mid \delta^{-1}(q_i^c,w_i)\subseteq S, |S| = n-k_0+1\right\}.
$$
So a set $S\in\mathscr{U}_{k_0-1,i}$ contains the pair of states $\delta^{-1}(q_i^c,w_i)$ and also $n-k_0-1$ of the remaining $n-2$ states in $Q$. Therefore
$$
|\mathscr{S}_{k_0,i}| = |\mathscr{U}_{k_0-1,i}| = \binom{n-2}{n-k_0-1},
$$
where the first equality is obtained by Lemma \ref{lemma:bijection}. Now note that also $\overline{\mathscr{S}}_{k_0,i}$ has the same size, which together with the inclusion $\mathscr{S}_{k_0,i} \subseteq \overline{\mathscr{S}}_{k_0,i}$ gives the equality $\mathscr{S}_{k_0,i} = \overline{\mathscr{S}}_{k_0,i}$ for all $i$.\\
\\
Now let $S$ be an arbitrary proper subset of $Q$ having size
$k_0$. Write $S$ as $\left\{q_{i_1},\ldots,q_{i_{k_0}}\right\}$.
Since $\sigma_W$ is a cyclic permutation on $Q$, there exists an
index $i_j$ such that $\sigma_W(q_{i_j})\notin S$, or equivalently
$q_{i_j}^c\notin S$. But then the complement $\tilde S :=
Q\setminus S$ satisfies
$$
q_{i_j}^c\in\tilde S,\qquad q_{i_j}\notin \tilde S\qquad\textrm{and}\qquad |\tilde S| = n-k_0.
$$
This means that $\tilde S\in\overline{\mathscr{S}}_{k_0,i_j} = {\mathscr{S}}_{k_0,i_j}$, which in turn implies that $\tilde S\in\mathscr{T}_{k_0}$. Hence $\mathscr{T}_{k_0}$ contains all subsets of $Q$ of size $n-k_0$, which means that $\overline{\mathscr{T}}_{k_0} \subseteq \mathscr{T}_{k_0}$, which gives the desired equality $\mathscr{T}_{k_0} = \overline{\mathscr{T}}_{k_0}$.\hfill $\Box$ \\
\\
Having collected these results, we are ready to provide proofs of our main theorems.\\
\\
\textbf{Proof of Theorem \ref{maintheorem1}.} Let $S$ be an arbitrary $k$-subset of $Q$. By Proposition
\ref{prop:maxsets}, the collection $\mathscr{T}_{n-k}$ contains all $k$-subsets of $Q$, so in particular $S\in\mathscr{T}_{n-k}$. By Lemma \ref{lemma:length}
there exists a word $w_{S}\in\Sigma^\star$ of
the form $w_{S} = w_{i_1}\ldots w_{i_{n-k}}$ for which
$\delta(Q,w_S) = S$.\hfill $\Box$ \\
\\
\textbf{Proof of Corollary 1.} Theorem \ref{maintheorem1} directly implies that every singleton can be reached by a word of the form $w_{i_1}\ldots w_{i_{n-1}}$. So in fact there exist at least $n$ synchronizing words of this form.\hfill $\Box$ \\
\\
\textbf{Proof of Theorem \ref{maintheorem2}.} The second statement is a direct consequence of Theorem \ref{maintheorem1}. Since the automaton is synchronizing, there exists at least one pair $\left\{q',q''\right\}\subseteq Q$ which can be mapped directly to a singleton. That is, for some label $a\in\Sigma$, we have $|\delta(\left\{q',q''\right\},a)|=1$. The pair $\left\{q',q''\right\}$ is reachable by a word $w$ of length at most $n(n-2)$, so the word $wa$ is synchronizing and has length at most $n(n-2)+1 = (n-1)^2$.\hfill $\Box$

\section{$1-$contracting automata and circularity}\label{section:circular}

In this section we make some comments on the scope of our results and on the relation between $1-$contractivity and circularity. An automaton is called circular if it has a label that cyclically permutes the states. 

%
%

As a special case of $1-$contracting automata we now consider
circular automata having a label $a\in\Sigma$ of deficiency $1$,
that is $|\delta(Q,a)| = |Q|-1$. Such automata are
$1-$contracting, since one can first read the $1-$deficient label
and then cyclically permute the excluded state. In this case, we
have the following result:

\begin{proposition} Let $\mathscr{A} = (Q,\Sigma,\delta)$ be a circular automaton with $n$ states. Suppose that $\mathscr{A}$ has a
$1$-deficient label $a$ and is circular with respect to $b$. Let $q_a$ be the state that is excluded
by $a$ and let $q_a^c$ be the unique contracting state for $a$.
Denote the distance on the circle from $q_a$ to $q_a^c$ by $d$. If
$\textrm{gcd}(d,n) = 1$, then
\begin{enumerate}
\item $\mathscr{A}$ is synchronizing.%
\item The shortest synchronizing word of $\mathscr{A}$ has length
at most $(n-1)^2$.%
\item For every nonempty set $S\subseteq Q$ of size $k$, there
exists a word $w_S$ of length at most $n(n-k)$ such that
$\delta(Q,w_S) = S$.
\end{enumerate}
\end{proposition}
\textbf{Proof.} Let $b$ be a label that cyclically permutes the
states. Let $q_1 = q_a$ and enumerate the other states according
to their position on the circle. Let $W = 
\left\{w_1,\ldots,w_n\right\}$
be defined by $w_i=ab^{i-1}, 1\leq
i\leq n$. Then $W$ forms an efficient $1-$contracting collection, since
$\delta(Q,w_i) = Q\setminus\left\{q_i\right\}$. Note that the label $b$ rotates the excluded and contracting state simultaneously, preserving the distance. Therefore, after reading the word $w_i$, the distance on the circle between the excluded and contracting state is equal to $d$. This immediately determines the state map $\sigma_W$:  
$$
\sigma_W(q) = \delta(q,b^d),\qquad \textrm{\ for all\ } q\in Q.
$$ 
Let $q\in Q$ and suppose $\textrm{gcd}(d,n) = 1$. Then $\sigma_W^n(q) = \delta(q,b^{nd}) = q$ and $\sigma_W^k(q) = \delta(q,b^{kd}) \neq q$ if $1\leq k <n$. In other words, $\sigma_W$ is a cyclic permutation on $Q$. 
Consequently the automaton is
aperiodically $1-$contracting and Theorem
\ref{maintheorem2} applies.\hfill $\Box$\\
\\
Pin \cite{pin2} already proved that synchronizing circular
automata with a prime number of states (later generalized by Dubuc
\cite{dubuc}) have a synchronizing word of length at most
$(n-1)^2$. Our result is slightly different, since we do not
assume that $\mathscr{A}$ is synchronizing. Furthermore, we get the additional property that all $k$-sets are
reachable with words of length at most $n(n-k)$.

In contrast with previous results on automata with some circular
structure, the class of $1-$contracting automata allows for
circles defined by a word consisting of different letters. To
illustrate this, we give an example.

\begin{example}
\emph{Consider the automaton $\mathscr{A}$ of Figure
\ref{fig:multletters}, where for simplicity we did not include the selfloops. In this automaton, the labels $a$ and $b$ are both permutations, and
together they make a circle through all states. With only these two labels the automaton would not be synchronizing. Therefore there is also a label $c$ which is $1$-deficient. Let $W$ be the collection of the words in the tables below. These words are $1-$deficient and their excluded and unique contracting states are indicated.
$$
\begin{array}{ccc}
\begin{array}{c|c|c}
\textrm{word} & \textrm{excl.} & \textrm{contr.} \\
\hline
w_1 = c & 1 & 2\\
w_2 = ca & 2 & 3\\
w_3 = cab & 3 & 4
\end{array}
&\qquad&
\begin{array}{c|c|c}
\textrm{word} & \textrm{excl.} & \textrm{contr.} \\
\hline
w_4 = cab^2 & 4 & 5\\
w_5 = cab^3 & 5 & 6\\
w_6 = cab^3a & 6 & 1
\end{array}
\end{array}
$$
$W$ is an efficient $1-$contracting collection. Clearly the state map $\sigma_W$ is a cyclic permutation, so it follows that $\mathscr{A}$ is aperiodically $1-$contracting. By Theorem \ref{maintheorem2}, all $k-$subsets of the state set are reachable by words of length at most $6(6-k)$. In particular, this automaton satisfies the \v{C}ern\'y conjecture.}\\
\phantom{asdf}\hfill $\blacksquare$
\end{example}

\begin{center}
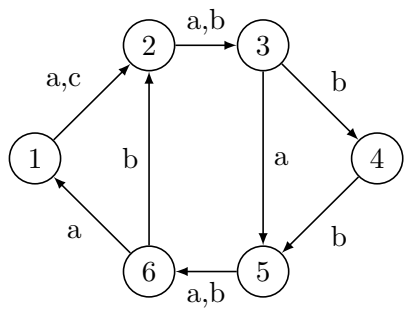

\begin{tikzpicture}[-latex ,auto ,node distance =1.5 cm and 1.5cm ,on grid ,
semithick , state/.style ={ circle ,top color =white , bottom
color = white!20 , draw, black , text=black , minimum width =.5
cm}]

\node[state] (2) {2};
\node[state] (3) [right =of 2] {3};
\node[state] (4) [below right=of 3] {4};
\node[state] (5) [below left =of 4] {5};
\node[state] (6) [left =of 5] {6};
\node[state] (1) [above left =of 6] {1};

\path (1) edge [bend right = 0] node[above left] {a,c} (2);
\path (2) edge node[above] {a,b} (3);
\path (3) edge node[right] {a} (5);
\path (5) edge node[below] {a,b} (6);
\path (6) edge node[below left] {a} (1);


\path (3) edge [bend right = 0] node[above right] {b} (4);
\path (4) edge [bend right = 0] node[below right] {b} (5);
\path (6) edge [bend right = 0] node[left] {b} (2);

\end{tikzpicture}
\captionof{figure}{An automaton where the labels $a$ and $b$ together make
a circle.}\label{fig:multletters}
\end{center}

Observe that the previous example can be easily extended to an infinite family of aperiodically $1-$contracting automata that have a circle with multiple labels. Theorem \ref{maintheorem2} applies as long as the excluded state and the contracting state move simultaneously along the circle on a mutual distance that is relative prime to the number of states, while the other states are permuting.

As a final illustration, we give an example that shows that our results are not restricted to automata with a cycle through all states. 
\begin{example} \emph{The automaton in Figure \ref{fig:withoutcircle} is aperiodically $1-$contracting with the following collection $W$ of $1-$deficient words:
$$
\begin{array}{ccc}
\begin{array}{c|c|c}
\textrm{word} & \textrm{excl.} & \textrm{contr.} \\
\hline
w_1 = a & 1 & 2\\
w_2 = ab & 2 & 4
\end{array}
&\qquad&
\begin{array}{c|c|c}
\textrm{word} & \textrm{excl.} & \textrm{contr.} \\
\hline
w_3 = abc & 3 & 1\\
w_4 = ac & 4 & 3
\end{array}
\end{array}
$$
Here the state map $\sigma_W$ is given by
$$
\sigma_W(1) = 2,\quad \sigma_W(2) = 4,\quad \sigma_W(3) = 1,\quad \sigma_W(4) = 3,
$$
which defines a cyclic permutation on $Q$. So also here Theorem \ref{maintheorem2} applies: all $k-$subsets
of the state set are reachable with a word of length at most
$4(4-k)$ and the \v{C}ern\'y conjecture holds true. In fact the
word $acabca$ is a synchronizing word that is obtained by
concatenating three of the $1-$deficient words $w_i$.}
\end{example}
\begin{center}
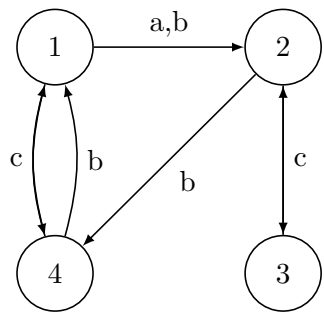

\begin{tikzpicture}[-latex ,auto ,node distance =3 cm and 3cm ,on grid ,
semithick ,
state/.style ={ circle ,top color =white , bottom color = white!20 ,
draw, black , text=black , minimum width =1 cm}]
\node[state] (1) {1};
\node[state] (2) [right =of 1] {2};
\node[state] (3) [below =of 2] {3};
\node[state] (4) [left  =of 3] {4};
\path (1) edge [bend right = 0] node[above] {a,b} (2);
\path (2) edge [bend right = 0] node[below right] {b} (4);
\path (4) edge [bend right = 15] node[right] {b} (1);
\path (2) edge [bend right = 0] node[right] {c} (3);
\path (3) edge (2);
\path (1) edge [bend right = 15] node[left] {c} (4);
\path (4) edge [bend left = 15]  (1);
\end{tikzpicture}
\captionof{figure}{An aperiodically $1-$contracting automaton without circle.}\label{fig:withoutcircle}
\end{center}

\section{Conclusion and perspectives}

We have introduced the notion of $1-$contracting automata in terms
of the existence of $1-$deficient words. We formulated an
aperiodicity property and have shown that this property suffices
to prove synchronization. If the $1-$deficient words with this
property have length at most $n$, then the automaton satisfies the
\v{C}ern\'y conjecture. In particular this is the case if the state map induced by a collection of minimal $1-$deficient words is a cyclic permutation. 

We believe that some elements in our definitions are not really
essential for the results, and that they only help in organizing a
clear proof. So a natural question is if the conditions can be
weakened, while still exploiting the ideas we presented in this
paper. We propose two directions in which this could be achieved:
\begin{enumerate}
\item Synchronization is a consequence of an automaton being aperiodically $1-$contracting. If
instead we assume synchronization, can we drop the aperiodicity assumption?

\item In order to prove the existence of a synchronizing word of length at most $(n-1)^2$, we require reachability of all subsets of size $n-1$ by words of length at
most $n$. Can this be relaxed by only requiring that any $(n-1)$-set has a subset which is reachable from $Q$ by a word of length at most $n$?
\end{enumerate}

Another question that naturally comes to mind is how to find a
synchronizing word in an $n-$state
aperiodically $1-$contracting automaton. From the proof
of Theorem \ref{maintheorem1}, it follows that there exists a synchronizing
word being the concatenation of $n-1$ of the $1-$deficient words.
However, it is not directly obvious how to formulate a full
constructive algorithm.

We proved for the class of aperiodically $1-$contracting automata that each reachable subset of size $k$ is reachable by a
word of length at most $n(n-k)$ if there exists an efficient $1-$contracting collection. In fact we believe that this is a
universal property, so we finish our paper with this conjecture:
\begin{conjecture}
Let $\mathscr{A} = (Q,\Sigma,\delta)$ be an $n$-state automaton.
If $S\subseteq Q$ is a set of size $k$ and there exists a word $w$
such that $\delta(Q,w) = S$, then there exists a word with this
property of length at most $n(n-k)$.
\end{conjecture}
This conjecture would imply the \v{C}ern\'y conjecture, since in a
synchronizing automaton there always exists a pair of states that
can be mapped to a singleton by some letter $a\in\Sigma$.

\end{document}